\documentclass[11pt]{amsart}

\usepackage{graphicx}
\usepackage{amsmath,amssymb,amsfonts}
\usepackage{amsthm}
\usepackage{mathrsfs}
\usepackage{tikz-cd}
\usepackage{hyperref}

\newcommand{\C}{\mathcal{C}}
\newcommand{\D}{\mathcal{D}}
\newcommand{\Real}{\mathbb{R}}
\newcommand{\OG}{\mathcal{O}_G}
\newcommand{\Csym}{\C^{\otimes}}
\newcommand{\GH}{\underline{G/H}}

\newcommand\opn{\operatorname}
\newcommand\op{\operatorname{op}}

\newcommand\Cat{\operatorname{Cat}}
\newcommand\Fun{\operatorname{Fun}}
\newcommand\THH{\operatorname{THH}}
\newcommand\THR{\operatorname{THR}}
\newcommand\colim{\operatorname{colim}}
\newcommand{\BOG}{\operatorname{BO}_n(G)}

\newcommand\bTop{\underline{\operatorname{\textsf{Top}}}}
\newcommand\Map{\operatorname{Map}}
\newcommand\bMap{\operatorname{\textsf{Map}}}
\newcommand\Mfld{\operatorname{Mfld}}
\newcommand\MG{\underline{\Mfld}^G}
\newcommand\MGB{\underline{\Mfld}^{G,f\opn{-fr}}}
\newcommand\Disk{\operatorname{Disk}}
\newcommand\DG{\underline{\Disk}^G}
\newcommand\DGB{\underline{\Disk}^{G,f\opn{-fr}}}
\newcommand\Spectra{\operatorname{\textsf{Spectra}}}
\newcommand\bSpectra{\underline{\operatorname{\textsf{Spectra}}}}

\theoremstyle{thmstyleone}
\newtheorem{theorem}{Theorem}[section]
\newtheorem{proposition}[theorem]{Proposition}

\theoremstyle{thmstyletwo}
\newtheorem{example}[theorem]{Example}

\theoremstyle{thmstylethree}
\newtheorem{definition}[theorem]{Definition}

\title{Applications of equivariant factorization homology}
\author{Aleksandar Miladinović}

\begin{document}

\begin{abstract}
In this paper we use the equivariant version of factorization homology constructed using the parametrized higher category theory and show that it can be used to describe the results used in the series of papers.
\end{abstract}

\maketitle
\tableofcontents

\section{Introduction}

Factorization homology represents homology theories for manifolds of fixed dimension.
It was introduced by Ayala and Francis in \cite{AF15} and Lurie in the form of chiral homology \cite{HA}.
It originated in the work of Beilinson and Drinfeld in \cite{BD04}, as well as in the work of Salvatore \cite{Sa01} and Segal \cite{Se10}.

There are numerous reasons to study factorization homology.
For one, factorization homology is used for defining topological quantum field theories \cite{CG16}.
Additionally, factorization homology allows us to define and study homology theories of manifolds in particular.

\subsection*{Higher category theory}
Factorization homology is a construction in higher category theory, or, to be more precise, in the setting of $\infty$-categories.

What are $\infty$-categories? To answer this question, we will ask a more general question:
What would be a correct category in which one could study homotopy theory?
Such a category should feature topological spaces as objects, continuous maps as morphisms,
 and should consider homotopies of continuous maps between topological spaces,
 as well as homotopies between homotopies, homotopies between homotopies of homotopies, etc.
It turns out that such category exists, if we accept to expand what we mean by category.
This has been explained in great detail in Lurie's book \cite{HTT}, hence we will jump straight to the definition of $\infty$-categories.

\begin{definition}
An $\infty$-category is a simplicial set $X$ such that the dotted lifts exist
\[
\begin{tikzcd} [row sep=2em, column sep=2em]
\Lambda^n_i \arrow[d,hookrightarrow]\arrow[r] & X \\
\Delta^n \arrow[ur,dashrightarrow] &
\end{tikzcd}
\]
for every $n\in\mathbb{N}$ and every $0<i<n$, where $\Delta^n$ is the standard $n$-simplex simplicial set, and where $\Lambda^n_i$ is the $i$-th horn of $\Delta^n$
i.e., $\Lambda^n_i$ is a subsimplicial set of $\Delta^n$ which features all $(n-1)$-faces of $\Delta^n$ containing $i$.
\end{definition}

Let $X$ be an $\infty$-category. The set of vertices $X_0$ is viewed as the set of objects while the set of edges $X_1$ is viewed as the set of morphisms.
From there the things get interesting.

Let $x$, $y$ and $z$ be three objects and let $f\colon x\to y$ and $g\colon y\to z$ be two morphisms in this $\infty$-category.
The diagram $x\xrightarrow{f} y\xrightarrow{g} z$ represents a horn $\Lambda^2_1$ in $X$.
By definition, such horn admits a lift i.e., we arrive at the diagram
\[
\begin{tikzcd} [row sep=2em, column sep=2em]
& y\arrow[rd,"g"] & \\
x \arrow[ur,"f"]\arrow[rr,"h"] & & z
\end{tikzcd}
\]
Given by the $2$-simplex $\sigma \colon  \Delta^2\to X$.
We can think of $\sigma$ as a homotopy form $g\circ f$ to $h$.
It is, however, important to note that there is no canonical choice of a composition arrow $g\circ f$, rather,
 the composition of $f$ and $g$ is defined \emph{up to homotopy}, which is part of the data.
We say that $\sigma$ exhibits $h$ as the composition of $f$ and $g$.
At first, this may seem like a disadvantage, however, it can be shown that the space of all possible compositions of $f$ and $g$ is in fact contractible,
 hence, the composition is unique up to homotopy.

To give more insight into the last statement let $\sigma'\colon \Delta^2\to X$ be the $2$-simplex of the form
\[
\begin{tikzcd} [row sep=2em, column sep=2em]
& y\arrow[rd,"g"] & \\
x \arrow[ur,"f"]\arrow[rr,"h'"] & & z
\end{tikzcd}
\]
Then, the maps $h$ and $h'$ will be homotopic in a suitable sense.
Let $\tau\colon \Delta^2\to X$ be the degenerate triangle given by
\[
\begin{tikzcd} [row sep=2em, column sep=2em]
& x\arrow[rd,"f"] & \\
x \arrow[ur,"id_x"]\arrow[rr,"f"] & & y
\end{tikzcd}
\]
We can further construct a $\Lambda^3_2$-horn in $X$ whose sides will be given by $\sigma$, $\sigma'$ and $\tau$
\[
\begin{tikzcd} [row sep=2em, column sep=2em]
& y\arrow[rrdd,"g"] & & \\
& x \arrow[u,"f"]\arrow[drr,"h'"'] & & \\
x \arrow[rrr,"h"']\arrow[ur,"id_x"']\arrow[uur,"f"] & & & z
\end{tikzcd}
\]
This horn can be filled to a $3$-simplex $\Delta^3\to X$.
In particular, we have the lower $2$-face of this simplex of the form
\[
\begin{tikzcd} [row sep=2em, column sep=2em]
& x\arrow[rd,"h'"] & \\
x \arrow[ur,"id_x"]\arrow[rr,"h"] & & h
\end{tikzcd}
\]
exhibiting a homotopy between $h$ and $h'$.

Furthermore, constructions such as limits and colimits in the $\infty$-category can be viewed as \emph{homotopy} limits and colimits in a suitable sense.

In addition, every category can be interpreted as an $\infty$-category via its \emph{nerve}.
For a small category $\mathcal{C}$, the nerve of $\C$ is a simplicial set $N(\C)$ given by
\[
N(\C)_n = \operatorname{Hom}_{\Cat}([n],\C)
\]
where $[n]$ is the category corresponding to the poset $\{1,2,\dots,n\}$ with its usual linear order.
Furthermore, if the morphism sets of a category $\C$ can be endowed with the structure of a topological space (or a simplicial set), one can make a construction of a so-called \emph{topological} (or \emph{simplicial}) nerve of a category $\C$, taking into account this structure.

Using this construction, we can build an $\infty$-category of spaces $\operatorname*{Spc}_{\infty}$
 as a topological nerve of the category of $\operatorname{CW}$-complexes and continuous maps between them.
Moreover, in such category, the notion of limits and colimits corresponds to the notion of \emph{homotopy} limits and colimits, which is exactly what we have wished for.

\subsection*{Parametrized higher category theory}
Let $X$ be a topological space with an action of a group $G$.
Let $\OG$ be the category of $G$-orbits i.e., the category of transitive $G$-spaces and $G$-equivariant maps between them.
We can think of elements of $\OG$ as cosets $G/H$ by choosing a basepoint of the orbit, with $H$ being the stabilizer of the basepoint.

The space $X$ can be represented by a functor
$F\colon \OG^{\op}\to \operatorname{Top}$, $G/H \mapsto X^H$,
hence, in order to study equivariant spaces categorically, one needs to consider $G$-spaces as functors depicted above.
The first problem that we encounter is that working with functors as objects is more complicated than working with (ordinary) categories.
Secondly, and more importantly, what happens when we want to incorporate homotopy into our framework?

In the same way that the higher category theory turned out to be the ideal setting in which one can do homotopy theory,
 parametrized higher category theory represents a good framework when dealing with the problems of equivariant homotopy theory.
The theory had been developed by Barwick and his students in a series of papers \cite{BDGNS16,Nar16,Nar17,NS22,Shah23}.
Some of the applications can be seen in \cite{Hor19,HHKWZ20}, as well as in the PhD thesis of the author \cite{Mil22}.

The general idea stems from the Grothendieck--Lurie correspondence:

\begin{theorem}
For a simplicial set $S$ there is an equivalence of $\infty$-categories
\[
u\colon  \Fun(S, \Cat_{\infty}) \xrightarrow{\simeq}\operatorname{Fib}^{\operatorname{coc}}(S)
\]
where $\Cat_{\infty}$ is the $\infty$-category of (small) $\infty$-categories, $\Fun(S, \Cat_{\infty})$ is the $\infty$-category
 of simplicial maps between $S$ and $\Cat_{\infty}$ and $\operatorname*{Fib}^{\operatorname{coc}}(S)$ is the $\infty$-category of coCartesian fibrations over $S$.
\end{theorem}

CoCartesian fibrations are special kind of morphisms, but since it will not be necessary to state its definition, it will be omitted.

Therefore, setting $S=\OG^{\op}$ we arrive to the definition of a $G$-$\infty$-category.
Namely, a $G$-$\infty$-category is a coCartesian fibration $\C\to\OG^{\op}$.
A functor of $G$-$\infty$-categories $\C\to\OG^{\op}$ and $\D\to\OG^{\op}$, or a $G$-functor,
 is a functor of $\infty$-categories $F\colon \C\to\D$ lying over $\OG^{\op}$, i.e., $F$ fits into a commutative diagram
\[
\begin{tikzcd} [row sep=2em, column sep=2em]
\C \arrow[rr,"F"]\arrow[dr] & & \D \arrow[dl] \\
& \OG^{\op} &
\end{tikzcd}
\]
and $F$ respects the coCartesian structure (in particular, $F$ sends coCartesian morphisms to coCartesian morphisms, see \cite{HTT} for more details).
Given a $G$-$\infty$-category $\C\to\OG^{\op}$, we can think of the fiber of the orbit $G/H$, denoted with $\C_{[G/H]}$, as the $\infty$-category of $H$-objects.
If $K\leq H$ are subgroups of $G$, then the lift of the map $G/K\to G/H$ gives us the functor $\C_{[G/H]}\to \C_{[G/K]}$
 which we can think of as a restriction functor.
Similarly, for an element $g\in G$ and $H\leq G$ we have a conjugation map $gHg^{-1}\to H$, which lifts to the functor $\C_{[G/H]}\to \C_{[G/gHg^{-1}]}$.

These examples illustrate why such a definition of an equivariant $\infty$-category is the right one: the $G$-$\infty$-category,
 in addition to the $G$-objects, contains $H$-object, for $H$ a subgroup of $G$, and all of them are linked in a suitable way (as we have seen).

There is also a notion of a $G$-symmetric monoidal $\infty$-category.
These categories are important since the $G$-factorization functor needs to be a $G$-symmetric monoidal functor
 in order to satisfy the equivalents of Eilenberg--Steenrod axioms for generalized homology theories.
Unfortunately, there are technical difficulties when working with $G$-symmetric monoidal categories.
In particular, $G$ cannot be any group.
Fortunately, $G$ can be a finite group, which is the most important case.
Additionally, $G$ is allowed to be a compact Lie group as long as the orbit category $\OG$ consists only of those orbit spaces with finite stabilizers.
Therefore, for the remainder of this paper, we will only consider $G$ to be a finite group or a compact Lie group,
 in which case $\OG$ is assumed to be an $\infty$-category of those transitive $G$-spaces with finite stabilizers.

\section{Equivariant factorization homology}

Equivariant factorization homology has two inputs: geometric and algebraic.
Geometric input is an equivariant manifold while the algebraic input is an equivariant disc algebra.
Let us start with the geometric input.

\subsection*{$\boldsymbol{G}$-Manifolds}
$G$-manifolds of a fixed dimension $n$ can be organized into a $G$-$\infty$-category in the following way:
consider a functor $\OG^{\op}\rightarrow \Cat_{\infty}$ sending the orbit $G/H$ to the $\infty$-category $N(\Mfld^H)$,
the nerve of the category of smooth $H$-manifolds and $H$-equivariant smooth embeddings between them.
Using the Grothendieck-Lurie correspondence we obtain the coCartesian fibration $\MG\rightarrow\OG^{\op}$ which we call the $G$-$\infty$-category of (smooth) $G$-manifolds.
Moreover, we can consider $G$-manifolds with additional tangential structure.
To be more precise, a smooth $G$-manifold $M$ is equipped with a tangent bundle map in the form $\tau \colon  M\rightarrow \BOG$.
If $B$ is a $G$-space and $f\colon B\to \BOG$ is a $G$-map, then the $B$-framing (or $f$-framing) on $M$ consists of a $G$-homotopy commutative diagram
\[
\begin{tikzcd} [row sep=2em, column sep=2em]
& B\arrow[rd,"f"] & \\
M \arrow[ur]\arrow[rr,"\tau"] & & \BOG
\end{tikzcd}
\]
In other words, the tangential structure on $M$ is determined by the map $f$.
Such $f$-framed $G$-manifolds can also be organized into a $G$-$\infty$-category $\MGB$.
Additionally, the disjoint union of $G$-manifolds provides a $G$-symmetric monoidal structure on $\MGB$.

\begin{example} \label{trivialframingex}
Let $B=\ast$ be just a point and let $f\colon \ast\to\BOG$ be a $G$-map.
A $f$-framing on a $G$-manifold $M$ corresponds to the trivialization of the tangent bundle of $M$.
To be more precise, tangent classifier map of $M$ factors through a point, hence, we obtain the following diagram
\[
\begin{tikzcd} [row sep=2em, column sep=2em]
\opn{TM}\arrow[d]\arrow[r] & V\arrow[d]\arrow[r] & \opn{EO}_n(G)\arrow[d] \\
M \arrow[r] & * \arrow[r] & \BOG
\end{tikzcd}
\]
where both inner rectangles are pullback diagrams, hence, the outer rectangle as well.
In such setting, $V$ is the $n$-dimensional $G$-representation, and since the left rectangle is a pullback diagram,
 the tangent bundle of $M$ is $\opn{TM}\cong M\times V$.
Instead of $f$-framed, we will often write $V$-framed in such case of framing on $G$-manifolds.
\end{example}

\subsection*{$\boldsymbol{G}$-Disk algebras}
$G$-disks represent the link between the algebra and the geometry of $G$-manifolds.
The $G$-$\infty$-category of $G$-disks is used for defining $G$-disk algebras
 which are again used as coefficients for equivariant version of factorization homology.
At the same time, $G$-disks provide insight into geometry of $G$-manifolds by capturing the local properties.
Furthermore, they can be linked with equivariant configuration spaces.
A $G$-disk can be viewed as a $G$-vector bundle $E\to U$, where $U$ is a finite $G$-set
(i.e., a coproduct of finite number of orbits).
We notice that $E$ is, in fact, a $G$-manifold, hence,
 we can consider a $G$-$\infty$-category of (framed) $G$-discs $\DGB$ as a full $G$-$\infty$-subcategory of $\MGB$ spanned by $G$-discs.
Furthermore, $\DGB$ admits a $G$-symmetric monoidal category inherited from $\MGB$.
Let $\Csym$ be a $G$-symmetric monoidal category.
An $f$-framed $G$-disc algebra is a $G$-symmetric monoidal functor $A\colon \DGB\to \Csym$.

\subsection*{$\boldsymbol{G}$-Factorization homology}
Since we have defined what are the input parameters, it is time to define $G$-factorization homology:

\begin{definition}
Let $M\in \MGB_{[G/H]}$ be an $H$-manifold (since it lays in the fiber of $G/H$),
 and let $A$ be an $f$-framed $G$-disk algebra taking values in a $G$-$\infty$-category $\C$.
The equivariant factorization homology of $M$ with coefficients in $A$, denoted with $\int_M A$, is given by the parametrized $\GH$-colimit
\[
\int_M A = \GH\text{-}\colim\Big(\DGB_{/\underline{M}}\to \GH\underline{\times}\DGB\to \GH\underline{\times}\C\Big)
\]
where $\GH = {\OG^{\op}}_{(G/H)/} \simeq \mathcal{O}_H^{\op}$, $\DGB_{/\underline{M}}$ is the parametrized slice category
 and where $\GH\underline{\times}\DGB$ is a category obtained via the pullback diagram
\[
\begin{tikzcd} [row sep=2em, column sep=2em]
\GH\underline{\times}\DGB \arrow[d]\arrow[r] & \DGB \arrow[d] \\
\GH \arrow[r] & \OG^{\op}
\end{tikzcd}
\]
(and similarly for $\GH\underline{\times}\C$).
\end{definition}

Parametrized colimits are fairly complicated and technically demanding, but the ideas behind them are natural \cite{Shah23}.
To be more precise, if we take $G$ to be a trivial group, then the $G$-parametrized colimits (or simply $G$-colimits) coincide with the regular colimits.

Alternatively, there is an adjunction of $G$-functors
\[
i_! \colon  \Fun_G(\DGB,\C)\rightleftarrows \Fun_G(\MGB,\C) \colon  i^*
\]
where $\Fun_G$ represents the $\infty$-category of $G$-functors between two $G$-$\infty$-categories,
 $i$ is the inclusion functor $i\colon \DGB\to\MGB$,
 $i^*$ is given by precomposition with $i$, and $i_!$ is given by the $G$-left Kan extension along $i$.
As it turns out  \cite[Section 4]{Hor19}, or \cite[Chapter 8]{Mil22}, the functor $i_!$ represents $G$-factorization homology functor.
Moreover, the $G$-factorization homology can be expanded to a $G$-symmetric monoidal functor.

\subsection*{Homology theory}

As stated in the introduction, $G$-factorization homology represents homology theory for $G$-manifolds
 and as such needs to satisfy some equivalents of Eilenberg--Steenrod axioms.
In particular, it needs to satisfy the $G$-$\otimes$-excision property and needs to respect $G$-sequential unions.
Let us break this down into more detail:

\begin{definition}\cite[5.1.1]{Hor19} or \cite[9.1.1]{Mil22}
Let $M$ be a $G$-manifold and let $[-1, 1]$ be a closed interval endowed with the trivial $G$-action.
By $G$-collar decomposition (or $G$-collar gluing), we mean a surjective equivariant map $f \colon  M \to [-1, 1]$
 such that the restriction $M|_{(-1,1)} \to (-1, 1)$ is a manifold bundle map with a choice of trivialization $M|_{(-1,1)}\cong M_0 \times (-1, 1)$ where $M_0 = f^{-1}(0)$.
We will denote with $M_+ = f^{-1}(-1,1]$ and $M_- = f^{-1}[-1,1)$.
\end{definition}

\begin{definition}
Let $F\colon  \MGB \to \C$ be a $G$-symmetric monoidal functor such that for every $G$-manifold $M$
 with a G-collar decomposition $f\colon  M \to [-1, 1]$ the induced map $F(M_-)\otimes_{F(M_0\times (-1,1))} F(M_+) \to F(M)$ is an equivalence,
 where $F(M_-)\otimes_{F(M_0\times (-1,1))} F(M_+)$ is given by the two-sided bar construction.
In such case we say that F satisfies the $G$-$\otimes$-excision property.
\end{definition}

\begin{definition}
Let $M$ be a $G$-manifold.
A $G$-sequential union of $M$ is a sequence of open $G$-submanifolds $M_1 \subseteq M_2 \subseteq\dots\subseteq M$ with $M = \bigcup_{i=1}^{\infty} M_i$.
\end{definition}

\begin{definition}
Let $F \colon  \MGB\to \C$ be a $G$-symmetric monoidal functor and $M =\bigcup_{i=1}^{\infty}M_i$ a $G$-sequential union of $M$.
Then $F$ induces a map
\[
\colim_i F(M_i)\to F(M).
\]
We say that $F$ respects $G$-sequential unions if this map is an equivalence.
\end{definition}

Finally, we have

\begin{proposition}\!\emph{\cite[9.2.4 and 9.3.3]{Mil22}}
Let $\C$ be a $G$-symmetric monoidal category and let $A\colon \DGB\to\C$ be an $f$-framed $G$-disk algebra.
Then the $G$-factorization functor $\int_- A \colon \MGB\to\C$ satisfies the $G$-$\otimes$-excision property and respects $G$-sequential unions.
\end{proposition}

Furthermore, the axiomatic characterization of $G$-factorization homology \cite[9.4.3]{Mil22} tells us that $G$-factorization homology accounts for all homology theories
 i.e., all $G$-symmetric monoidal functors that satisfy the $G$-$\otimes$-property and respect $G$-sequential unions.

\subsection*{Naive $\boldsymbol{G}$-factorization homology}
If one wishes to disregard the parametri\-zed structure and consider only the $\infty$-category of $f$-framed $G$-manifolds,
 which would be the $\infty$-category $\MGB_{[G/G]}$ when $G$ is finite, then the construction of the factorization homology is still possible \cite{Wee20}.
In such construction there is no \emph{genuine} action of $G$
i.e., we have no restriction and conjugation functors since there is no parametrized structure.
Furthermore, a functor from such category does not carry any information about the $G$-action, it is just a functor from a category with $G$-objects to some other category.
Hence, we call this construction the \emph{naive} $G$-factorization homology.
Even though it seems like a step down from the parametrized point of view, this construction does provide some interesting examples, as we shall see.

\begin{example} \cite[3.12]{AF15}
The first example that we are going to see is the simplest one,
 hence we have left it out of the section dedicated to the applications of the equivariant factorization homology.
Let $G$ be a trivial group.
In this case the construction of the $G$-factorization homology agrees with the one provided by Ayala and Francis in \cite{AF15}.
Furthermore, let $M$ be the unit interval $[-1,1]$.
Note that $M$ is a manifold with a boundary, hence, we can consider the $\infty$-category of one-dimensional oriented disks (possibly) with boundary $\Disk^{\opn{or}}_1$.
Let $\Csym$ be a symmetric monoidal category which is $\otimes$-presentable and $A\colon \Disk^{\opn{or}}_1\to \Csym$ a one-dimensional oriented disc algebra.
Let us denote (by abuse of notation) with $A = A((-1,1))$, and with $L = A([-1,1))$, $R = A((-1,1])$.
Embeddings $(-1,1)\hookrightarrow [-1,1)$ and $(-1,1)\hookrightarrow (-1,1]$ give $R$ and $L$ a structure of a right and left $A$-module, respectively.
Furthermore, $M$ admits obvious collar gluing
\[
[-1,1]\cong [-1,1) \cup_{(-1,1)} (-1,1].
\]
Therefore, the excision property of factorization homology gives $\int_{[-1,1]} A \simeq L \otimes_A R$.
Readers not familiar with the notion of a $\otimes$-presentable $\infty$-category can safely ignore it,
 since it is a technical term, and does not stop us from understanding the idea behind this example.
\end{example}

\section{Applications}

In the final section of this paper we will look at some (and certainly not all) examples of applications of equivariant factorization homology in the literature.

\subsection{Ordinary factorization homology}

As stated in the previous example, when $G$ is taken to be a trivial group, the equivariant factorization homology can be viewed as the (ordinary) factorization homology.
There are many applications, as depicted in \cite{AF19}, hence, we will not dive into every one of them, but will select a few.

\subsubsection{Ordinary and generalized homology theories}
Let $B=*$ be a point and let $f\colon \ast\to \opn{BO}_n$, let $\Spectra$ be the $\infty$-category of spectra,
 and let $A\colon \Disk^{f\text{-fr}}_n\to \Spectra$ be an $f$-framed $n$-disk algebra.
By \ref{trivialframingex} the $f$-framing on a manifold $M$ corresponds to the trivialization of the tangent bundle $\opn{TM}\cong M\times \mathbb{R}^n \simeq M$,
 hence, the factorization homology gives $\int_M A \simeq \Sigma^{\infty} M \wedge A(\mathbb{R}^n)$
where $\Sigma^{\infty} M$ is the suspension spectrum of $M$.
In the special case when $A(\mathbb{R}^n)$ is chosen to be the Eilenberg--MacLane spectrum $H\mathbb{Z}$ we obtain a spectrum which represents the ordinary homology theory $H_{\bullet}(M)$.
Moreover, if $C$ is any Abelian group and $HC$ the Eilenberg--MacLane spectrum of $C$,
 factorization homology produces a spectrum which corresponds to the ordinary homology of $M$ with coefficients in $C$, $H_{\bullet}(M;C)$.

\subsubsection{Hochschild homology} \label{hochschildhomology}
\begin{figure}[h]
\begin{tikzpicture}
\begin{scope}[thick, every edge/.style = {draw,->}, every node/.style={sloped,allow upside down}]
\draw[red] (0,0) to [out=90, in=180] (2,2);
\draw[red] (0,0) to [out=-90, in=180] (2,-2);
\draw[red] (2,2) to [out=0, in=90] (3.9,0.2);
\draw[red] (2,-2) to [out=0, in=-90] (3.9,-0.2);

\draw[blue] (3.7,0) to [out=90, in=0] (2,1.7);
\draw[blue] (3.7,0) to [out=-90, in=0] (2,-1.7);
\draw[blue] (2,1.7) to [out=180, in=90] (0.4,0.2);
\draw[blue] (2,-1.7) to [out=180, in=-90] (0.4,-0.2);

\draw[green] (0.6, 0.2) to [out=90, in=180] (2,1.4);
\draw[green] (2,1.4) to [out=0, in=90] (3.4,0.2);

\draw[green] (0.6, -0.2) to [out=-90, in=180] (2,-1.4);
\draw[green] (2,-1.4) to [out=0, in=-90] (3.4,-0.2);
\end{scope}
\end{tikzpicture}
\caption{Collar gluing of the circle}
\label{slika}
\end{figure}
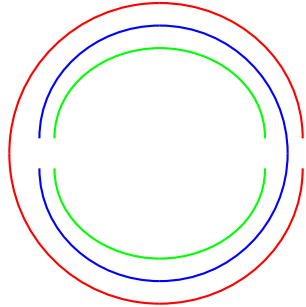

Another great example of factorization homology is the connection with Hochschild homology (see \cite[3.19]{AF15} or \cite[3.31]{AF19}).
Let $M$ be a circle $S^1$.
If we write $S^1$ as the set $\{(x,y)\in \Real^2 \mid x^2+y^2=1 \}$,
 then the projection on the $x$-axis provides us with a collar gluing of $S^1\cong \Real \cup_{\Real \sqcup \Real} \Real$ presented as in Figure \ref{slika}.
Let $\Csym$ be a nice enough symmetric monoidal $\infty$-category (in particular, $\Csym$ needs to be $\otimes$-presentable $\infty$-category),
 and let $A\colon \Disk^{\opn{or}}_1 \to \Csym$ be an oriented $1$-disk algebra.
Note that the orientation on the category of one-dimensional manifolds is given by a suitable choice of framing \cite[2.2]{AF15}.
In literature, $A$ is also called an associative algebra object in $\Csym$.
We obtain
\[
\int_{S^1} A \simeq A(\Real) \otimes_{A(\Real)\otimes A(\Real)^{\op}} A(\Real)
\]
where the superscript ``op" stems from the opposite orientation of the second $\Real$ in $\Real \sqcup \Real$ from the standard orientation on Euclidean spaces.
Additionally, this is equivalent to the Hochschild homology associated to the algebra $A$
i.e., $A(\Real) \otimes_{A(\Real)\otimes A(\Real)^{\op}} A(\Real) \simeq \opn{HH}_{\bullet}(A)$.
Furthermore, we can make the connection between the factorization homology and the topological Hochschild homology ($\THH$) by taking $A$ to be an associative ring spectrum.
The action on the circle translates to the action on $\THH(A)$ which serves as a motivation for studying factorization homology in the equivariant setting,
 keeping in mind that there is an action of the group $O(2)$ on $S^1$ and therefore on $\THH(A)$.

\subsection{Bredon homology and Borel equivariant homology}

As ordinary factorization homology can encode general homology theories via spectra, similar can be said about the equivariant factorization homology.
Let $G$ be a finite group.
In \cite[Chapter 5]{Bre67}, Bredon defined equivariant (co)homology theories and later proved that,
 for every functor $A\colon \OG\to\opn{Ab}$ (where $\opn{Ab}$ is the category of Abelian groups)
 there is a unique equivariant (Bredon) homology theory $H^G_{\bullet}(-;A)$ such that $H^G_0(G/H;A)=A(G/H)$.
Each such theory is constructed as homology of a certain chain complex $C^G_{\bullet}(-;A)$.
By \cite[Section 5.1]{Wee20} every Bredon homology theory can be realized via equivariant factorization homology
\[
H^G_{\bullet}(M;A) \simeq \int_M C^G_{\bullet}(-;A)
\]
where $C^G_{\bullet}(-;A)\colon  \DG\to\opn{Ch}^{\oplus}_{\mathbb{K}}$ is a $G$-disk algebra
 with coefficients in the category of chain complexes over a field $\mathbb{K}$.
Moreover, for every $G$-disk $E$ viewed as the total space of a vector bundle $E\to G/H$, we have $E\simeq G/H$ and hence $H^G_0(E;A)=A(G/H)$.
In particular, if $A$ is chosen to be a constant coefficient system at $\mathbb{Z}$
i.e., $A(G/H)=\mathbb{Z}$ for every $G/H\in\OG$, then the Bredon homology coincides with the notion of Borel equivariant homology.

\subsection{Real topological Hochschild homology}

Let $C_2$ be the cyclic group of order $2$, and let $\sigma$ be its one-dimensional sign representation (i.e., $\sigma$ is $\Real$ with an action of $C_2$, sending $x$ to $-x$).
Furthermore, regard $S^1$ represented as the set $\{(x,y)\in\Real^2 \mid x^2+y^2=1\}$ with $C_2$ action as the reflection with respect to the $y$-axis.
Then there is a $C_2$-equivariant collar gluing on $S^1$ given by $S^1\cong \sigma\bigcup_{\Real\sqcup_{C_2}\Real}\sigma$.
Let $\bSpectra^{C_2}$ be the $G$-symmetric monoidal category of genuine $C_2$-spectra
 and let $\underline{\Disk}^{C_2,\sigma\text{-fr}}$ be the $C_2$-symmetric monoidal category of $\sigma$-framed one-dimensional disks.
Then, for $A\colon \underline{\Disk}^{C_2,\sigma\text{-fr}}\to \bSpectra^{C_2}$ we have the result $\int_{S^1} A \simeq \THR(A)$ of Horev \cite[7.1.1,\,7.1.2]{Hor19}
where $\THR(A)$ is the topological real Hochschild homology of $A$.
Moreover, $\THR(A)$ admits an action of $C_2$ which can be refined to $O(2)$ action by \cite[10.3]{Mil22}.

\subsection{Twisted topological Hochschild homology}
Let $C_n$ be a cyclic group of order $n$, and let $\underline{\Disk}^{C_n,\Real\text{-fr}}$ be the $C_n$-symmetric monoidal category
 of $\Real$-framed one-dimensional $C_n$-disks, with $C_n$ acting trivially on $\Real$.
Let $S^1$ be a circle with standard action of $C_n$ (by rotations) and let $A\colon \underline{\Disk}^{C_n,\Real\text{-fr}}\to \bSpectra^{C_n}$.
By abuse of notation, let us write $A:=A(\Real)$.
Then, by \cite[7.2.3]{Hor19}, we have
\[
\bigg(\int_{S^1} A \bigg)^{\Phi C_n} \simeq \THH(A;A^{\tau})
\]
where the superscript $\Phi C_n$ denotes the geometric fixed points,
 $A^{\tau}$ is the $A-A$-bimodule with twisted multiplication given by first acting on the scalar by the generator $\tau\in C_n$:
\[
A\otimes A^{\tau} \otimes A \to A,\;\;
x\otimes a\otimes y \mapsto \tau x \otimes a\otimes y
\]
and where $\THH(A;A^{\tau})$ is the topological Hochschild homology of $A$ with coefficients in $A^{\tau}$.
In, particular, this result tells us that $\THH(A;A^{\tau})$ admits a natural circle action.

\subsection{Equivariant nonabelian Poincaré duality}
Equivariant factorization homology can be used to prove the equivariant version of nonabelian Poincaré duality:
Let $G$ be a finite group, $V$ a $G$-representation and $M$ a $V$-framed $G$-manifold.
Additionally, let $X$ be a pointed $G$-space such that $\pi_k(X^H)=0$ for all subgroups $H<G$ and all $k< \dim(V^H)$.
It can be shown \cite{HHKWZ20} that the $G$-functor $\bMap_*((-)^+, X)\colon \underline{\Mfld}^{G,V\text{-fr}} \to \bTop_*^{G}$ is a homology theory of $G$-manifolds
(i.e., it satisfies the $G$-$\otimes$-excision property and respects $G$-sequential unions)
(where the superscript $+$ denotes the one-point compactification, and where $\bTop_*^{G}$ represents the $G$-$\infty$-category of pointed $G$-spaces).
Hence, we obtain a natural equivalence \cite[Theorem 2.2]{HHKWZ20}
\[
\int_M \Omega^V X \simeq \Map_*((M)^+, X)
\]
where we can write $\bMap_*((-)^+, X)\colon \underline{\Disk}^{G,V\text{-fr}} \to \bTop_*^{G}$
 for a $V$-framed $G$-disk algebra with $\bMap_*((V)^+, X)_{[G/G]}\simeq \Map_*((V)^+, X) = \Omega^V X$.
In the case when $G$ is a trivial group, the upper equivalence can be upgraded into
\[
\int_M \Omega^n X \simeq \Map_c(M,X)
\]
where $M$ is a manifold of dimension $n$, $X$ is $n-1$-connective space and $\Map_c(M,X)$ is the space of compactly supported maps from $M$ to $X$.
Furthermore, when $n=1$, this equivalence reduces to the Goodwillie's quasi-isomorphism $\opn{HH}_{\bullet}(\Omega X) \simeq LY$,
 where $LX = \Map(S^1,X)$ is the free loop space of $X$ (see \ref{hochschildhomology} and \cite{Good85}).

\subsection{Norm construction}
The norm construction is the vital tool for understanding the multiplicative structure on the genuine equivariant stable category when $G$ is finite.
It was first introduced by Hill, Hopkins and Ravenel in their paper on the Kervaire invariant problem \cite{HHR16}.
When $G$ is a compact Lie group the situation is different and more complicated.
Even though there have been some constructions of norms of particular compact Lie groups ($S^1$ and $O(2)$ to be precise)
 in \cite{ABGHL18} and \cite{AGH21}, up until recently, there has not been a general construction of the norm for compact Lie groups.
The first such construction is proposed in \cite{BHM22} using factorization homology.
Let $G$ be a compact Lie group of dimension $n$ and let $U$ be a $G$-universe i.e., a countable sum of some set of irreducible $G$-representations.
Furthermore, let $A$ be a spectrum.
By abuse of notation, let us write $A\colon \Disk^{*\text{-fr}}_n\to \Spectra$ for a point-framed $n$-disk algebra with $A(\Real^n)=A$.
We define the absolute norm of $A$, $N^G_e A$ as
\[
N^G_e A = I^U_{U^G}\int_G A
\]
where $I^U_{U^G}$ denotes the point-set change of universe functor from $\Real^{\infty}= U^G$ to~$U$,
 and where $\int_G A$ can be viewed as a $G$-object by the natural left $G$-action on $G$.
This construction gives us a genuine $G$-spectrum object $N^G_e A$ indexed on the universe~$U$, starting from a spectrum object $A$.
Moreover, this construction can be upgraded to a relative norm construction $N^G_H$ for any closed subgroup $H<G$.
In the case when $G=S^1$, this construction agrees with the relative norm constructed in \cite{ABGHL18}.

\bibliographystyle{amsplain}

\end{document}